\documentclass[12pt]{amsart}
\usepackage{amsmath,amssymb,amsthm}
\newtheorem*{theorem}{Theorem}

\newcommand{\SOCL}{sequence of constituent lengths}
\newcommand{\AFS}{\mathrm{AFS}}
\newcommand{\B}{\mathrm{B}}
\newcommand{\F}{\mathbb{F}}
\newcommand{\N}{\mathbb{N}}
\newcommand{\gr}[2]{{#1}^{#2}}

\begin{document}
\mbox{}
\title{(Finite) presentations of Bi-Zassenhaus loop algebras}
\author{G. Jurman}
\address{Fondazione Bruno Kessler, via Sommarive 18, I-38123 Povo (Trento), Italy}
\email{jurman@fbk.eu}
\thanks{The author is a member of UMI and INdAM-GNSAGA}

\begin{abstract}
We prove that Bi-Zassenhaus loop algebras are finitely presented up to
central and second central elements. In particular, we show an explicit finite
presentation for a Lie algebra whose quotient over its second centre is 
isomorphic to a Bi-Zassenhaus loop algebra. 
\end{abstract}

\maketitle

\section{Introduction}
\label{sec:intro}
The family of simple modular Lie algebras nowadays known as
Hamiltonian algebras was discovered in the fifties by Albert
and Frank \cite{AlbFra54} and it was named after them. By definition, 
they are graded over a non-cyclic elementary abelian
$p$-group. Shalev \cite{Sha94b} noticed that they also have a
cyclic grading with a non-singular derivation cycling over this grading. 
This property allowed him to build an infinite-dimensional loop algebra starting from the original simple one. 
The new structure is of maximal class, \textit{i.e.} it is $\N$-graded with all homogeneous
components of dimension one except the first one which has dimension
two and generates the entire algebra. 
These new algebras are known as Albert-Frank-Shalev ($\AFS$, for
short) algebras.
In the following years there has been a growing interest in these
algebras. Caranti, Mattarei and Newman developed new
techniques to construct more algebras \cite{CarMatNew97} 
and later Caranti and Newman achieved a classification theorem for odd
characteristic fields of definition \cite{CarNew00}. The
remarkable result they reached yields that every infinite-dimensional
$\N$-graded Lie algebras of maximal class generated by its first
homogeneous component can be built starting from an $\AFS$-algebra.
To get the claim they need a theorem proved by Carrara in \cite{Carr98,Carr01}, which states that
every $\AFS$ algebra is uniquely determined by a suitable finite
dimensional quotient. This result was obtained as a corollary of a more
general property she discovered: $\AFS$-algebras are
finitely presented up to central and second central elements.
A similar classification for characteristic two is still valid \cite{Jur05}, 
but in this case another family of algebras is
involved, the Bi-Zassenhaus loop algebras ($\B_l$ for short).
The loop construction process for them is exactly the same as for Albert-Frank-Shalev
algebras, but the simple finite-dimensional algebra at the beginning is different. 
Here we show that the property described above for
$\AFS$-algebras is satisfied by $\B_l$-algebras, too. 
The fact that $\B_l$-algebras are uniquely determined by a suitable
finite-dimensional quotient again follows as a corollary.
In particular, the proven result is the following
\begin{theorem}
For every Bi-Zassenhaus loop algebra $\B_l (g,h)$, there
exists a finitely presented graded Lie algebra $M(g,h)$ such that $M(g,h) / Z_2 (M(g,h)) \cong \B_l(g,h)\ .$
\end{theorem}
The construction of the algebra $M$ by means of cohomological arguments is shown in~\cite{Jur04}.
Since the centre of $M$ is infinite-dimensional, a group-theoretical result of
B.H.~Neumann~\cite [pp. 52-53]{Rob93} transferred to Lie algebra shows
that this implies that $\B_l$ itself is \emph{not} finitely presented.
The relations in the presentation of the algebra $M$ are
retrieved by expanding several suitable generalized Jacobi identities,
thus machine computations took a main role throughout the work. 
These computations were performed by the software $p$-Quotient 
Program~\cite{HavNewObr97} developed at the Australian National University in Canberra.
Nevertheless, all the proofs are independent from such calculations.
\section{Preliminaries}
\label{sec:preliminaries}
A graded Lie algebra 
\begin{displaymath}
L= \bigoplus _{i=1}^{\infty} L_i\ , 
\end{displaymath}
defined over a field of positive characteristic $p$ and generated by
$L_1$ is said to be of \emph{maximal class} 
if $\dim (L_1 )=2$ and $\dim (L_i )\leq 1$ for $i>1$. 
Its elements are written left-normed and in exponential form:
\begin{displaymath}
[u\gr{v}{n}] =  [ [ [u\underbrace{v]v]\cdots v}_{n}]\ .
\end{displaymath}
The \emph{two-step centralizers} $C_i$ of $L$ are defined as the
one-dimensional subspaces $C_{L_1}(L_i)$ of $L_1$
centralizing the homogeneous components $L_i$ when $i>1$, while $C_1 =
C_2$ is formally assumed. Following the classical notation, we define
the element $y$ by means of $C_2 = \F y$ and we choose another 
element $x\in L_1\setminus \F y$. The pair $\{x,y\}$ is thus a set of generators for $L_1$. 
Let $C_k$ the first occurrence of a two-step centralizer different from $\F y$. 
When the class of $L$ is larger than $k+1$, all the two-step
centralizers coincide with $C_2$ apart from isolated occurrences of
different subspaces, so it is possible to define a \emph{constituent}
of $L$ as a subsequence $(C_i,\dots,C_{j})$ where all centralizers
coincide with $C_2$ but the last one, and either $i=1$ or
$C_{i-1}\not=C_2$. Homogeneous elements lying in class $i-1$ and $j$
are respectively said to be \emph{at the beginning} and \emph{at the end} of the
constituent and they are not centralized by $C_2$.

A graded Lie algebras of maximal class is determined by its 
sequence of two-step centralizers, and, when only two of them are
distinct, by its \SOCL . The first constituent of a graded Lie
algebra of maximal class whose class is larger than $k+2$ is always 
of length $2q$ for some power $q=p^h$ called its \emph{parameter}, and all other constituents
can only be \emph{short} (\textit{i.e.} of length $q$) or \emph{long}
($2q$) or \emph{intermediate} ($2q-p^s$, where $0\leq s\leq h-1$). 
The standard result~\cite [Prop. 5.6]{CarMatNew97} holds also when central and second central elements occur, as shown in~\cite{Carr01}: for a (large enough) algebra of parameter $q$,  
\begin{equation}
\textsl{the only possible constituent lengths are $2q$ and $2q-2^s,\; 0\leq s\leq h$\ .} \tag{$\mathcal{CL}$}
\end{equation}
From now on, unless explicitly stated, we will assume $p=2$. This
allows us to ignore signs safely so that the generalized Jacobi
identity reads as follows in the two equivalent forms:
\begin{displaymath}
 [v  [w \gr{u}{\lambda}]]  = \sum_{i=0}^{\lambda}  {\binom{\lambda}{i}}  [v \gr{u}{i} w\gr{u}{\lambda - i}] = \sum_{i=0}^{\lambda}  {\binom{\lambda}{i}}  [v \gr{u}{\lambda - i} w \gr{u}{i}]\ .
\end{displaymath}
The generalized Jacobi identity can be more effectively used by introducing the element $z=x+y$: 
if $0 \ne v \in L_{i}$, for some $i>1$, and if 
$C_{i},C_{i+1}, \dots, C_{i+n-1}\in \{\F x,\F y\},$
then for every non-zero commutator $[v x_{1} x_{2} \dots x_{n}]$ with $x_{i} \in \{ x, y \}$ we have that
$[v x_{1} x_{2} \dots x_{n}] = [v \gr{z}{n}] \ .$

The Bi-Zassenhaus loop algebras $\B_l(g,h)$ form a family of elementary objects among infinite-dimensional graded Lie algebras of maximal class. 
They are defined when the characteristic of the underlying field is two and $(h,g-1)\in\N\times\N$. 
If exponential notation is employed in order to indicate consecutive
occurrences of patterns, the \SOCL\ of the algebra
$\B_l(g,h)$ reads as
\begin{displaymath}
2q,2q-1,\left( (2q)^{\eta-1}, (2q-1)^2  \right)^\infty\ ,
\end{displaymath}
where $q=2^h$ and $\eta=2^g-1$.
Thus in a $\B_l$-algebra the only intermediate constituents are those
of maximal length $2q-1$ and there are no short constituents. 
Suppose $\displaystyle M=\bigoplus _{i=1}^{\infty }
M_i $ is a graded Lie algebra such that $\displaystyle M/Z_2 (M) $ is a
graded Lie algebra of maximal class. We define a {\it constituent} of $M$
as a constituent of $ M/Z_2 (M) $, ignoring the central or second central elements of $M$. 

The notation $w^{-c}$ will be used for an homogeneous element such that $[w^{-c}\gr{x}{c}] = w\ .$
Then we define the elements
\begin{displaymath}
v_n=
[y\gr{x}{2q-1}(y\gr{x}{2q-2}(y\gr{x}{2q-1})^{\eta-1}y
\gr{x}{2q-2})^n]\ ,        
\end{displaymath}
which lie in the homogeneous components of indices $2q+dn$ where
$d=2^{g+h+1}-2$ is the dimension of the simple algebra $\B(g,h)$
and $n$ ranges over the non-negative integers. When the two parameters
$a,b$ lie in the intervals $2\leq a\leq h+1$ and $h+2\leq b\leq g+h$,
we can define the following elements, consistently with the notation in~\cite{Jur04}:
\begin{displaymath}
\begin{array}{lclr}
\theta_n^1 &=&  [v_n x ] & (2q+1+dn)\phantom{\ .} \\
\theta_n^a &=&  [v_n y\gr{x}{2q -2^{h+2-a}-1}y] & (4q-2^{h+2-a}+1+dn)\phantom{\ .}\\
\theta_n^b &=&  [v_n y\gr{x}{2q-2} (y\gr{x}{2q-1})^{\eta-2^{g+h+1-b}}y\gr{x}{2q-2}y] & (2q(\eta+3-2^{g+h+1-b})-1+dn)\phantom{\ .} \\
\theta_n^\omega &=&  [\theta^1_{2n+1}y] & (2q+2+d(2n+1))\ .
\end{array}
\end{displaymath}
The number in parentheses is the \emph{weight} of
the element, \textit{i.e.} the index $i\geq 1$ of the homogeneous component
$M_i$ in which it lies. Finally, we define the shorthands
\begin{eqnarray*}
\mu_{n,i} &=& 
\left\{
\begin{array}{ll}
{} [v_n y\gr{x}{2q-3}] & \textrm{for $i=1$}\ , \\
{} [v_n y\gr{x}{2q-2}(y\gr{x}{2q-1})^{i-2}y\gr{x}{2q-2}] & \textrm{otherwise}\ .
\end{array} 
\right. 
\end{eqnarray*}
Note that $[\mu_{0,\eta-2^{g+h+1+b}+2}y]= \theta_0^b$ when $h+2\leq b\leq g+h$.
\section{The finite presentation}
\label{sec:finite_presentation}
Fix a particular $\B_l$ algebra $L$ and then find a finite set
$R'$ of homogeneous relations such that $M'=\langle x ,y : { R'} \rangle $ is a presentation of 
a suitable $m$-dimensional Lie algebra such that $M'/Z_2(M')$ is 
isomorphic to a graded quotient of $L$ and remove the relations of degree $m+1$. 
What obtained will be shown to be a finitely presented, infinite-dimensional
graded Lie algebra $M$ such that $M/Z_2 (M)\cong L $. Adding to $R'$
the set of relations used to annihilate the central and second
central elements, a graded Lie algebra of maximal class $M$ gets defined,  
isomorphic to $L$ and uniquely determined by a suitable
finite-dimensional quotient $M'$.
When all the homogeneous relations that define the algebra $L$ up to central and
second central elements up to class 
$m=2q(\eta+2)$ are added to $R'$, then the resulting algebra $M$ starts as a graded Lie algebra of maximal
class with initial segment of constituent lengths 
$2q, 2q-1, 2q^{\eta-1},2q-1$
and $C_i \in \{ \F x, \F y \}$ for every $1\leq i\leq m$, but most of the $m-2$ relations that in every class $2\leq i\leq m$ set the two-step centralizer $C_i$ are actually redundant.
Thus the proof of the main theorem reduces to show that the following $q+h+\eta$ relations
\begin{displaymath}
\left\{
\begin{array}{ll}
[y\gr{x}{2j+1}y]=0, &\qquad 0\leq j \leq q-2\ ,  \\
 {} [\theta_0^t x]=0, &\qquad 1\leq t\leq g+h\ \textrm{or}\  t=\omega \ , \\
 {} [\mu_{0,t+2} y]=0,
&\qquad 0\leq t\leq \eta-2\;\textrm{and $\eta-t\not =
   2^\alpha$, for $1\leq \alpha \leq g-1$}  \ .
\end{array}
\right.
\end{displaymath}
defining the set $R$ are sufficient to give a presentation for $M(g,h)$.

The first set of $q$ relations defines the parameter $q=2^h$:
\begin{displaymath}
 [y\gr{x}{2j+1}y]=0 \quad\textrm{for $0\leq j\leq q-2$}\ ,\qquad
 [y\gr{x}{2q+1}]=0\ . 
\end{displaymath}
This makes all homogeneous components have the same two-step
centralizer up to weight
$2q$. Moreover, they force all further constituents to
be short, long or intermediate. The next relation
\begin{displaymath}
 [y\gr{x}{2q-1}y\gr{x}{q-1}yx]=0
\end{displaymath}
states that the second constituent is not short: in particular, this
implies that no more short constituents or two-step centralizers other
than the first two will be involved. 
A standard argument shows that the second constituent
cannot be long, so it can be only intermediate. The following $h-1$ 
relations 
\begin{displaymath}
 [y\gr{x}{2q-1}y\gr{x}{2q-2^s-1}yx]=0 \quad\textrm{for $1\leq s\leq h-1$}
\end{displaymath}
establish its length as the maximal possible $2q-1$. So the algebra is
not inflated and starts moving on the branch of the limit algebra
$\AFS(h,h+1,\infty,2)$: furthermore, we know that
from now on its \SOCL\ will contain only long and maximal intermediate
constituents. By the classification we know that the number of long constituents before two intermediate ones is a power of two
minus one or minus two. 
The remaining $\eta$ relations determine the algebra as belonging to the $\B_l$ family rather than to the $\AFS$ one:
\begin{gather*}
[y\gr{x}{2q-1}y\gr{x}{2q-2}(y\gr{x}{2q-1})^t y\gr{x}{2q-2}y(x)]=0 
\quad\textrm{for $0\leq t\leq \eta-2$}\ , \\
[y\gr{x}{2q-1}y\gr{x}{2q-2}(y\gr{x}{2q-1})^{\eta-1} y\gr{x}{2q-1}yx]=0\ , 
\end{gather*}
where the $(x)$ appears only when $\eta-t$ is a power of two.
\section{The proof: the structure}
\label{sec:proof_structure}
We will prove that the \SOCL\ of $(M,R)$ and $L$ are the same 
and that the two algebras differ only in central or second central
elements; in particular, we claim that the elements $\theta_n^t$ are
the only central and second central elements that may occur: the
construction shown in~\cite{Jur04} completes the picture.
For clarity, the expansion of the following Jacobi identities
which are not immediate is carried out in Appendix~\ref{sec:proof_expansions}.
\subsection{Using the first group of relations}
The first homogeneous component is obviously two-dimensional:
\begin{displaymath}
M_1 = \langle x,y \rangle\ .
\end{displaymath}
The first group of relations written above shows the structure of $M$
up to class $2q+1$:
\begin{displaymath}
M_i = \langle  [y\gr{x}{i-1}] \rangle\quad\textrm{for}\  2\leq i\leq
2q\ , \qquad 
M_{2q+1}=  \langle  [y\gr{x}{2q-1}y],  [y\gr{x}{2q}]=\theta_0^1 \rangle\ .
\end{displaymath}
In fact, the elements
\begin{displaymath}
 [y\gr{x}{j}y] \quad\textrm{for $0\leq j\leq 2q-2$}
\end{displaymath}
vanish either by the relations in $R$ (in the case $j$ odd) or by the following inductive argument (when $j$ is even):
\begin{displaymath}
0 =  [ [y\gr{x}{\frac{j}{2}}] [y\gr{x}{\frac{j}{2}}]] =
 {\binom{\frac{j}{2}}{\frac{j}{2}}} [y\gr{x}{j}y]\ .
\end{displaymath}
Finally, the relation in class $2q+2$ 
\begin{displaymath}
 [y\gr{x}{2q+1}]= [\theta_0^1 x]= 0
\end{displaymath}
and the expansions
\begin{eqnarray*}
0&=&  [y\gr{x}{2q-2} [yxy]]= [y\gr{x}{2q-1}yy]\ , \\
0 &=&  [ [y\gr{x}{q}] [y\gr{x}{q}]] =  {\binom{q}{q-1}}  [y\gr{x}{2q-1}yx] +
{\binom{q}{q}}  [y\gr{x}{2q}y] =  [\theta_0^1 y]
\end{eqnarray*}
show that the first constituent of $M$ has length $2q$, that
$C_{2q}=\F x$ and that $\theta_0^1$ is central:
\begin{displaymath}
M_{2q+2} = \langle  [y\gr{x}{2q-1}yx] \rangle\ .
\end{displaymath}

\subsection{The constituent lengths}
The relations just obtained are those needed to prove that property $(\mathcal{CL})$ holds in $M$.  As a direct consequence, if $v$ is not centralized by $y$, then
\begin{displaymath}
 [v_n y\gr{x}{k}y]=0\quad \textrm{when $k\not = 2q-2^s-1$ for $0\leq s\leq h$}\ . 
\end{displaymath}
By using this equation, it is not difficult to show that the second
constituent of $M$ is of maximal intermediate length $2q-1$, that its
last element is centralized by $x$ and that the elements $\theta_0^a$
are central, for $2\leq a\leq h+1$:
\begin{eqnarray*}
M_{2q+i} &=& 
\left\{
\begin{array}{ll}
\parbox[t]{4cm}{
$\langle  [y\gr{x}{2q-1}y\gr{x}{i-1}]\ ,$\\  
$\phantom{\langle}[y\gr{x}{2q-1}y\gr{x}{i-2}y]=\theta_0^a\rangle$
}
& 
\parbox[t]{6cm}{
for $3\leq i\leq 2q-1\ ,$ \\
$i = 2q-2^{h+2-a}+1$,\ $2\leq a \leq h+1\ ,$
}\\
\langle  [y\gr{x}{2q-1}y\gr{x}{i-1}]
\rangle & \textrm{otherwise}\ ,
\end{array}
\right.
\end{eqnarray*}

\begin{displaymath} 
M_{4q} =  \langle  [y\gr{x}{2q-1}y\gr{x}{2q-2}y] \rangle\ .
\end{displaymath}
Since two consecutive values of $i$ cannot be both of the type
$2q-2^{h+2-a}-1$, the thesis follows from the fact that every element
\begin{displaymath}
 [y\gr{x}{2q-1}y\gr{x}{j}y]
\end{displaymath}
vanishes by the property ($\mathcal{CL}$) when $0\leq j\leq 2q-3$ and
$j\not =2q-2^s-1$ for $1\leq s\leq h$ and otherwise, by the relation
\begin{displaymath}
 [y\gr{x}{2q-1}y\gr{x}{j}yx] =  [\theta_0^{h+2-s} x] = 0
\end{displaymath}
and the expansions
\begin{displaymath}
0 =  [y\gr{x}{2q-1}y\gr{x}{j-1} [yxy]] =
 [y\gr{x}{2q-1}y\gr{x}{j}yy]=   [\theta_0^{h+2-s} y]\ ,
\end{displaymath}
\begin{displaymath}
0 =  [ [y\gr{x}{2q-1}] [y\gr{x}{2q-1}]] = {\binom{2q-1}{0}} [y\gr{x}{2q-1}y\gr{x}{2q-1}]\ .
\end{displaymath}

\subsection{The number of long constituents} 
Now we deal with the remaining structure of $M$ up to class $m=2q(\eta+2)$, by showing that 
when $i$ runs between zero and $\eta-1$ and $k$ between zero and $2q-1$, with $k\not
= 2q-1$ for $i=\eta-1$, we have that 
\begin{displaymath}
M_{4q+2qi+k}=
\left\{
\begin{array}{ll}
\parbox[t]{6.1cm}{
$\langle  [y\gr{x}{2q-1}y\gr{x}{2q-2}(y\gr{x}{2q-1})^{i}y\gr{x}{2q-1}]\
,$\\
$\phantom{\langle}[y\gr{x}{2q-1}y\gr{x}{2q-2}(y\gr{x}{2q-1})^i y\gr{x}{2q-2}y]=\theta_0^b \rangle$
\\
} &
\parbox[t]{5cm}{
for $k=2q-1$\ ,\\
$i=\eta-2^{g+h+1-b}$, \textrm{and}\\
$h+2\leq b\leq g+h$\ ,
}\\
\langle  [y\gr{x}{2q-1}y\gr{x}{2q-2}(y\gr{x}{2q-1})^i y\gr{x}{k}] \rangle &
\textrm{otherwise}\ .
\end{array}
\right.
\end{displaymath}
To prove the above result, we consider separately some cases:
\begin{itemize}
\item when $k=0$ 
\begin{itemize}
\item and $i>0$ (the case $i=0$ has been already dealt with in the previous subsection), first of all we have the relation $ [\theta_0^1 x]=[y\gr{x}{2q+1}]=0$:
\begin{displaymath}
0 =  [y\gr{x}{2q-1}y\gr{x}{2q-2}(y\gr{x}{2q-1})^{i-1}y\gr{x}{2q-2} [\theta_0^1 x]] =  {\binom{2q+1}{1}} [y\gr{x}{2q-1}y\gr{x}{2q-2}(y\gr{x}{2q-1})^i y\gr{x}{2q-1} x]\ ;
\end{displaymath}
\item in particular, if $i=\eta-2^{g+h+1-b}$, for $h+2\leq b\leq g+h$, then
the previous homogeneous component has dimension two; in this case, to show that $\theta_0^b$ is central, we can use the relation
\begin{displaymath}
 [\theta_0^b x]=0
\end{displaymath}
and the standard expansion
\begin{displaymath}
0 =
 [y\gr{x}{2q-1}y\gr{x}{2q-2}(y\gr{x}{2q-1})^{i-1}y\gr{x}{2q-3} [yxy]]=
 [\theta_0^b y]\ ;
\end{displaymath}
\end{itemize}
\item when $k=2q-2^s$, for $1\leq s\leq h$,
\begin{itemize}
\item and $i=0$, we can use the following Jacobi expansion: 
\begin{equation}
\label{first}
0 =  [ [y\gr{x}{2q-1}y\gr{x}{q-2^{s-1}-1}] [y\gr{x}{2q-1}y\gr{x}{q-2^{s-1}-1}]] =  [y\gr{x}{2q-1}y\gr{x}{2q-2} y\gr{x}{k-1}y]\ ;
\end{equation}
\item and $i>0$, by using the relation $ [\theta_0^{s+1} x]=0$ in the
  following identity, which is the case $n=0$ of ~\eqref{long2}:
\begin{displaymath}
0 =  [\mu_{0,i}  [\theta_0^{s+1}x]] =  [y\gr{x}{2q-1}y\gr{x}{2q-2}(y\gr{x}{2q-1})^i y\gr{x}{k-1}y]\ .
\end{displaymath}
\end{itemize}
\item {when $k=2q-1$ and $i\leq \eta-2$, the thesis explicitly follows by the relation 
$[\mu_{0,i+2} y]=0$, when $i\not = \eta-2^{g+h+1-b}$,
  while otherwise there is nothing to prove since the homogeneous
  component is two-dimensional.}
\item by proposition ($\mathcal{CL}$) in the remaining cases.
\end{itemize}
So we have proved that $M$ has the same structure of $L$ up to class
$m$, apart from some central elements.

\subsection{From the quotient algebra to the whole algebra}
Now we have to prove that the finite-dimensional quotient we built
determines $M$ as a $\B_l$-algebra, apart from the central and second
central elements $\theta_n^t$, for any integer $n\geq 1$.

The last equation in the previous section shows that
\begin{displaymath}
M_{2q+d} = \langle v_1 \rangle\ ;
\end{displaymath}
by induction the homogeneous component in class $2q+dn$ is
one-dimensional too, generated by the element
\begin{displaymath}
M_{2q+dn} = \langle v_n \rangle\ .
\end{displaymath}
In what follows we describe the structure of the next $d$ homogeneous
components of $M$, \textit{i.e.} of an entire period of the
algebra. 
The homogeneous component in class $2q+1+dn$ is two-dimensional:
\begin{displaymath}
M_{2q+1+dn}= \langle  [v_n x]=\theta_n^1,  [v_n y] \rangle \ ,
\end{displaymath}
where $\theta_n^1$ can be central or second central, since in the next
class 
\begin{displaymath}
M_{2q+2+dn} = 
\left\{
\begin{array}{ll}
\langle  [v_n yx] \rangle & \textrm{when $n$ is even}\ , \\
\langle  [v_n yx], [v_n xy]=\theta_{\frac{n-1}{2}}^\omega \rangle & \textrm{when $n$ is odd}\ .
\end{array}
\right.
\end{displaymath}
In fact, in addition to the standard identities
\begin{displaymath}
0 =  [v^{-1}_n  [yxy]] =  [v_n yy]
\end{displaymath}
and the following one derived by the relation $[\theta_0^1 x]=[y\gr{x}{2q+1}]=0$
\begin{displaymath}
\begin{split}
0 &= [v_{n-1}y\gr{x}{2q-2}(y\gr{x}{2q-1})^{\eta-2}y\gr{x}{2q-2} [\theta_0^1 x]] \\
&=  {\binom{2q+1}{1}} [v_{n-1}y\gr{x}{2q-2}(y\gr{x}{2q-1})^{\eta-1}y\gr{x}{2q}] = [\theta_n^1 x]=  [v_n xx]\ ,
\end{split}
\end{displaymath}
when $n=2s$ is even we have the equation
\begin{displaymath}
 [v_n xy] =  [\theta_n^1 y] = 0\ ,
\end{displaymath}
given by the expansion
\begin{equation}
\label{n_even}
 [ [v^{-(q-1)}_s] [v^{-(q-1)}_s]]=0
\ ,
\end{equation}
that shows the centrality of $\theta_1^1$ in those cases.

Now we consider class $2q+3+dn$: if $n\geq 1$ is odd, first we have the equation
\begin{displaymath}
0 =  [v_n  [yxy]] =  [\theta_{\frac{n-1}{2}}^\omega y]=  [v_n xyy]\ ;
\end{displaymath}
then, if $n=1$ we have the relation
\begin{displaymath}
0 =  [\theta_0^\omega x] =  [v_n xyx]\ ,
\end{displaymath}
otherwise we can use the expansion
\begin{equation}
\label{second}
0 =  [v_{n-2} y\gr{x}{2q-2}(y\gr{x}{2q-1})^{\eta-2}y\gr{x}{2q-2}
 [\theta_0^\omega x]] =  [\theta_{\frac{n-1}{2}}^\omega x] =  [v_n xyx]\ ,
\end{equation}
to show that the element $\theta_{\frac{n-1}{2}}^\omega$ is central.
Finally, if $q=2$ the element $ [v_n yxy]$ is just $\theta_n^2$, while when 
$q>2$ the relation $ [yxxxy]=0$ holds, so we can expand
\begin{displaymath}
0 =  [v^{-2}_n  [yxxxy]] =  {\binom{3}{2}}  [v_n yxy]  =  [v_n yxy]\ .
\end{displaymath}

Summarizing, we have
\begin{displaymath}
M_{2q+3+dn} = 
\left\{
\begin{array}{ll}
\langle  [v_n yxx], [v_n yxy]=\theta_n^2 \rangle &\textrm{when $q=2$}\ ,\\ 
\langle  [v_n yxx] \rangle & \textrm{otherwise}\ .
\end{array}
\right.
\end{displaymath}

Now we focus the attention on the homogeneous components up to class
$4q+dn$. 

If $q>2$, then $2q+3+dn < 4q-1 +dn$ so we can study what happens
in classes $2q+1+k+dn$, where $3\leq k\leq 2q-3$. We claim that
\begin{displaymath}
M_{2q+1+k+dn} =
\left\{
\begin{array}{ll}
\langle  [v_n y\gr{x}{k}], [v_n y\gr{x}{k-1}y]=\theta_n^a \rangle & \textrm{when $k=2q-2^{h+2-a}$}\\
&\textrm{for some $2\leq a\leq h+1$\ ,}\\
\langle  [v_n y\gr{x}{k}] \rangle & \textrm{otherwise}\ .
\end{array}
\right.
\end{displaymath}
If $k=2q-2^{h+2-a}$ for some $2\leq a\leq h+1$, by hypothesis 
$M_{2q+1+k+dn}$ is generated by both $ [v_n y\gr{x}{k}]$ and
$ [v_n y\gr{x}{k-1}y]$: since $h+2-a>0$ the number $k+1$ cannot be of the
form $2q-2^{h+2-a'}$ for any $2\leq a'\leq h+1$, so $ [v_n y\gr{x}{k}y]=0$
by ($\mathcal{CL}$); moreover, the following equation~\eqref{theta^ax}
\begin{equation}
\label{theta^ax}
0 =  [v_{n-1}y\gr{x}{2q-2}(y\gr{x}{2q-1})^{\eta-2}y\gr{x}{2q-2} 
 [\theta_0^{a} x]] =  [\theta_n^{a}x]\ ,
\end{equation}
and the standard expansion
\begin{displaymath}
0 =  [v_n y\gr{x}{2q-2^{h+2-a}-2} [yxy]] =  [\theta_n^{a}y]
\end{displaymath}
show that the element $ [v_n y\gr{x}{k-1}y]=\theta_n^a$ is central. 

The remaining case is covered by the proposition ($\mathcal{CL}$).

Now look at the class $4q-1+dn$, which in the case $q=2$
(\textit{i.e. $h=1$}) coincide with the class $2q+3+dn$:
\begin{displaymath}
M_{4q-1+dn} = \langle  [v_n y\gr{x}{2q-2}],
 [v_n y\gr{x}{2q-3}y]=\theta_n^{h+1}\rangle \ .
\end{displaymath}
The centrality of the element $\theta_n^{h+1}$ has been already proven; moreover, 
the expansion
\begin{equation}
\label{Xi}
0 =  [\Xi_n \Xi_n] =  [v_n y\gr{x}{2q-1}]
\end{equation}
where
\begin{eqnarray*}
\Xi_n &=& 
\left\{
\begin{array}{ll}
v_s & \textrm{for $n=2s$}\ , \\
{} [v_s y\gr{x}{2q-2}(y\gr{x}{2q-1})^\frac{\eta-1}{2}]
& \textrm{for $n=2s+1$}\ ,
\end{array}
\right.
\end{eqnarray*}
shows that
\begin{displaymath}
M_{4q+dn} = \langle  [v_n y\gr{x}{2q-2}y]\rangle = \langle [v_n y\gr{x}{2q-2}y]\rangle \ .
\end{displaymath}

Finally, the last homogeneous components are generated as follows, where $i$ runs
between zero and $\eta-1$ and $k$ between zero and $2q-1$, with $k\not
= 2q-2$ for $i=\eta-1$: 
\begin{displaymath}
M_{4q+2qi+k+dn}=
\left\{
\begin{array}{ll}
\parbox[t]{6cm}{
$\langle  [v_n y\gr{x}{2q-2}(y\gr{x}{2q-1})^{i}y\gr{x}{2q-1}]\ ,$\\
$\phantom{\langle}  [v_n y\gr{x}{2q-2}(y\gr{x}{2q-1})^i y\gr{x}{2q-2}y]=\theta_n^b \rangle$ 
\\
}
&
\parbox[t]{5cm}{
for $k=2q-1$\ , \\
$i=\eta-2^{g+h+1-b}$,\\ 
$h+2\leq b\leq g+h$\ ,
}\\
\langle  [v_n y\gr{x}{2q-2}(y\gr{x}{2q-1})^i y\gr{x}{k}] \rangle &
\textrm{otherwise}\ .
\end{array}
\right.
\end{displaymath}

To prove this result, we distinguish some cases.
\begin{itemize}
\item when $k=0$ and $1\leq i\leq \eta-1$ (the case $i=0$ has been already proven), then first we have that
\begin{displaymath}
0 =  [v_n y\gr{x}{2q-2} (y\gr{x}{2q-1})^{i} \theta_0^1] =  {\binom{2q}{0}} [v_n y\gr{x}{2q-2} (y\gr{x}{2q-1})^{i} y\gr{x}{2q-1} x]\ .
\end{displaymath}
In particular, in the case $i=\eta-2^{g+h+1-b}$ for some integer $h+2\leq b\leq
g+h$, two more identities are required, since $M_{4q+2qi+dn-1}$ is
two-dimensional: the former reads as 
\begin{displaymath}
0 =  [v_n y\gr{x}{2q-2}(y\gr{x}{2q-1})^{i} y\gr{x}{2q-3} [yxy]] =  [v_n y\gr{x}{2q-2}(y\gr{x}{2q-1})^{i} y\gr{x}{2q-2}yy]= [\theta_n^b y]\ ,
\end{displaymath}
while the latter uses the relation $ [\theta_0^b x]=0$:
\begin{equation}
\label{pag15}
0 =  [v_{n-1}y\gr{x}{2q-2}(y\gr{x}{2q-1})^{\eta-2} y\gr{x}{2q-2} [\theta_0^b x]]=  [v_n y\gr{x}{2q-2}(y\gr{x}{2q-1})^{i} y\gr{x}{2q-2}yx] =  [\theta_n^b x]\ ;
\end{equation}
\item when $k=2q-2^s$ for $1\leq s\leq h$, 
\begin{itemize}
\item and $i=0$, then the relation 
$[\mu_{0,2} y]=0$ can be used:
\begin{equation}
\label{pag16}
0 = [v_{n-1}y\gr{x}{2q-2}(y\gr{x}{2q-1})^{\eta-2}y\gr{x}{2q-2^s-1} [\mu_{0,2} y] ] = [v_n y\gr{x}{2q-2}y\gr{x}{k-1}y]\ ;
\end{equation}
\item and $i\geq 1$, we can use the relation $ [\theta_0^{s+1} x]=0$:
\begin{equation}
\label{long2}
0 = [\mu_{n,i})  [\theta_0^{s+1} x]] =  [v_n y\gr{x}{2q-2}(y\gr{x}{2q-1})^i y\gr{x}{k-1}y]\ ,
\end{equation}
\end{itemize}
\item when $k=2q-1$, 
\begin{itemize}
\item and $i=\eta-2^{g+h+1-b}$ for some $h+2\leq b\leq g+h$, then
  $M_{4q+2qi+2q-1+dn}$ is two dimensional;
\item and $i$ is not one of the above values, then define $\lambda$ 
as the exponent of the highest power of two dividing $i+1$ and use
the relation $[\mu_{0,i+2^\lambda+1} y]=0$ as follows:
\begin{equation}
\label{lambda}
0 = [v_{n-1}y\gr{x}{2q-2}(y\gr{x}{2q-1})^{\eta-1-2^\lambda}y\gr{x}{2q-2} [\mu_{0,i+2^\lambda+1} y] ] =  [v_n y\gr{x}{2q-2}(y\gr{x}{2q-1})^i y\gr{x}{2q-2}y]\ ;
\end{equation}
\end{itemize}
\item when $k$ is not one of the above values, then proposition ($\mathcal{CL}$) proves the claim.
\end{itemize}
Now to conclude look at the homogeneous component occurring for $k=2q-2$ and $i=\eta-1$:
\begin{displaymath}
M_{4q+2q(\eta-1)+2q-2+dn} = M_{2q+d(n+1)}= \langle [v_n y\gr{x}{2q-2}(y\gr{x}{2q-1})^{\eta-1}y\gr{x}{2q-2}] \rangle = \langle [v_{n+1}] \rangle\ . 
\end{displaymath}
\section*{Acknowledgement}
The author is grateful for his help to A.~Caranti, advisor of the doctoral dissertation \cite{Jur98} this work is based on.
He is also grateful to M.F.~Newman and S.~Mattarei for making useful suggestions and reading various versions of this paper and to an anonymous referee for her/his precious help in indicating improvements and signaling mistakes. 
\bibliographystyle{amsalpha}
\bibliography{bibliography}
\appendix
\section{The expansions}
\label{sec:proof_expansions}
\addtocounter{subsection}{-1}
\subsection{Computing tools}
\label{subsec:computing_tools}
Lucas' Theorem \cite{KnuWil89,Luc78} will be used several times: if $\displaystyle{a=\sum_{i=0}^n a_i\cdot 2^i}$
and $\displaystyle{b=\sum_{i=0}^n b_i\cdot 2^i}$ are the $2$-adic
expansions of two integers, then
\begin{displaymath}
{\binom{a}{b}}\equiv \prod_{i=0}^n {\binom{a_i}{b_i}} \pmod 2 \ .
\end{displaymath}
For instance, as a consequence of the formula above, we have that
\begin{displaymath}
{\binom{2^w-1}{i}}\equiv 1\pmod 2\quad\forall\; 0\leq i\leq 2^w-1\ .
\end{displaymath}
The following identity ($\mathcal{I}$) will be be useful in simplifying some evaluations: if $Q$ is a power of $2$, then 
\begin{equation}
\sum_{j=0}^s {\binom{(Q-1)s+r}{(Q-1)j+k}} \equiv \binom{r}{k} \pmod 2 \quad 0\leq r,k\leq Q-2\ . \tag{$\mathcal{I}$}
\end{equation}
To prove the identity ($\mathcal{I})$, consider for $z>0$ the power sum of the elements of the field $\mathbb{F}_{Q}$:
\begin{displaymath}
\sum_{\alpha\in\mathbb{F}_Q^*} \alpha^z \equiv
\begin{cases}
1 \pmod 2 & \textrm{if $Q-1 | z$,} \\
0 \pmod 2 & \textrm{otherwise\ .}
\end{cases}
\end{displaymath}
Then, for $n=(Q-1)s+r>0$ we have that
\begin{displaymath}
\sum_{\alpha\in\mathbb{F}_Q^*} (1+\alpha)^n \alpha^{-k}=
\sum_{\alpha\in\mathbb{F}_Q^*} \sum_{i=0}^n \binom{n}{i} \alpha^{i-k} =
\sum_{i=0}^n \binom{n}{i} \sum_{\alpha\in\mathbb{F}_Q^*} \alpha^{i-k} =
\sum_{j=0}^{s} \binom{n}{(Q-1)j+k}\ .
\end{displaymath}
Since $\alpha^n = \alpha^r$ in $\mathbb{F}_Q$, the leftmost term in the above equation only depends on $r=n\pmod{Q-1}$: thus, this is true for the rightmost term too, and then identity ($\mathcal{I}$) follows.
A more general statement, generalization of a congruence originally shown by Glaisher in 1899, is proven as Prop. 6 in \cite{Mat06} by multisection of series.

\subsection{Expansion of eq. \eqref{first}}
The first equation we deal with is~\eqref{first}, whose expansion gives
\begin{eqnarray*}
0 &=& [ [y\gr{x}{2q-1}y\gr{x}{q-2^{s-1}-1}] [y\gr{x}{2q-1}y\gr{x}{q-2^{s-1}-1}]] \\
&=&  [ [y\gr{x}{2q-1}y\gr{x}{q-2^{s-1}-1}] [y\gr{z}{3q-2^{s-1}-1}]]\\
&=&  [y\gr{x}{2q-1}y\gr{x}{q-2^{s-1}-1}\gr{z}{3q-2^{s-1}-1}y] +  [y\gr{x}{2q-1}y\gr{x}{q-2^{s-1}-1}\gr{z}{3q-2^{s-1}}] \cdot {\binom{3q-2^{s-1}-1}{q+2^{s-1}-1}} \\
&=&  [y\gr{x}{2q-1}y\gr{x}{q-2^{s-1}-1}\gr{z}{3q-2^{s-1}-1}y] \\
&=&  [y\gr{x}{2q-1}y\gr{x}{2q-2} y\gr{x}{k}y]\ ,
\end{eqnarray*}
since the binomial coefficient is equivalent to
\begin{displaymath}
{\binom{3q-2^{s-1}-1}{q+2^{s-1}-1}} \equiv 
{\binom{2}{1}}{\binom{q-2^{s-1}-1}{2^{s-1}-1}}\equiv 0\pmod 2\ .
\end{displaymath}

\subsection{Expansion of eq. \eqref{n_even}} 
Expanding from the back:
\begin{eqnarray*}
0 &=& [ [v_{s-1}
y\gr{x}{2q-2}(y\gr{x}{2q-1})^{\eta-1}y\gr{x}{q-1}]
 [v_{s-1}
y\gr{x}{2q-2}(y\gr{x}{2q-1})^{\eta-1}y\gr{x}{q-1}]] \\
&=& [ [v_{s-1}
y\gr{x}{2q-2}(y\gr{x}{2q-1})^{\eta-1}y\gr{x}{q-1}]  [y\gr{z}{ds+q}]]\\
&=& [v_{s-1}
y\gr{x}{2q-2}(y\gr{x}{2q-1})^{\eta-1}y\gr{x}{q-1}\gr{z}{ds+q}y]  \\
&\phantom{=}&+  [v_{s-1}
y\gr{x}{2q-2}(y\gr{x}{2q-1})^{\eta-1}y\gr{x}{q-1}\gr{z}{ds+q+1}]\cdot
\left( 
\sum_{i=1}^{\eta}  \binom{ds+q}{2qi} \right. \\
&&\qquad +  \binom{ds+q}{2q(\eta+1)-1} +  \sum_{j=0}^{s-2} \left(
 \binom{ds+q}{dj-2+2q(\eta+2)}  \right. \\
&&\qquad + \left. \left. \sum_{i=1}^{\eta-1} \binom{(2^{g+h+1}-2)s+2^{h}}{dj-2+2q(\eta+2+i)} +
 {\binom{(2^{g+h+1}-2)s+2^{h}}{dj-3+2q(2\eta+2)}}\right)\right)\\
&=&  [\theta_n^1 y]\ .
\end{eqnarray*}
In fact, the second and the last term in the above coefficient are equivalent
to zero modulo two in view of Lucas' Theorem, since their denominator is odd while the numerator
is even. The remaining terms can be rewritten in terms of $g$ and $h$ as follows:
\begin{eqnarray*}
\binom{ds+q}{2qi} &=&
\binom{(2^{g+h+1}-2)s+2^{h}}{2^{h+1}i} \\
\binom{ds+q}{dj-2+2q(\eta+2)} &=&
\binom{ds+q}{(2^{g+h+1}-2)(j+1)+2^{h+1}} \\
\binom{ds+q}{dj-2+2q(\eta+2+i)}
&=&
\binom{(2^{g+h+1}-2)s+2^{h}}{(2^{g+h+1}-2)(j+1)+2^{h+1}(i+1)}\ ,
\end{eqnarray*}
and thus arranged in an unique sum that can be shown to vanish by using Lucas' Theorem and identity ($\mathcal{I}$):
\begin{displaymath}
\sum_{j=0}^{s-1}
\sum_{i=1}^{2^g-1} 
\binom{(2^{g+h+1}-2)s+2^{h}}{(2^{g+h+1}-2)j+2^{h+1}i} \equiv 
\sum_{j=0}^{s-1}
\sum_{i=1}^{2^g-1} 
\binom{(2^{g+h}-1)s+2^{h-1}}{(2^{g+h}-1)j+2^{h}i}\equiv 0 \pmod 2\ ,
\end{displaymath}
since $2^{h-1}<2^h i$ for $i$ ranging between $1$ and $2^g-1$.

\subsection{Expansion of eq. \eqref{second}}
Then we have the equation~\eqref{second}, needed to prove the
centrality of $\theta_{\frac{n-1}{2}}^\omega$: note that, in this
case, the relation can be conveniently written as 
\begin{displaymath}
 [\theta_0^\omega x] =  [y\gr{x}{2q-1}y\gr{x}{2q-2}(y\gr{x}{2q-1})^{\eta-1}y\gr{x}{2q-1}yx ]
=  [y\gr{z}{2q(\eta+1)+2q-2}yx]\ ,
\end{displaymath}
since the element $ [y\gr{z}{2q(\eta+1)+2q-3}yxx]$ is non zero. The expansion begins as:
\begin{eqnarray*}
0 &=&  [v_{n-2} y\gr{x}{2q-2}(y\gr{x}{2q-1})^{\eta-2}y\gr{x}{2q-2}
 [\theta_0^\omega x]] \\
&=&  [v_{n-2} y\gr{x}{2q-2}(y\gr{x}{2q-1})^{\eta-2}y\gr{x}{2q-2} 
 [y\gr{x}{2q-1}y\gr{x}{2q-2}(y\gr{x}{2q-1})^{\eta-1}y\gr{x}{2q-1}yx ]]\\
&=&  [v_{n-2} y\gr{x}{2q-2}(y\gr{x}{2q-1})^{\eta-2}y\gr{x}{2q-2}
 [y\gr{z}{2q(\eta+1)+2q-2}yx]]\\
&=&  [v_{n-2} y\gr{x}{2q-2}(y\gr{x}{2q-1})^{\eta-2}y\gr{x}{2q-2}
 [y\gr{z}{2q(\eta+1)+2q-2}y]x]\\
&\phantom{=}&+  [v_{n-2} y\gr{x}{2q-2}(y\gr{x}{2q-1})^{\eta-2}y\gr{x}{2q-2}x
 [y\gr{z}{2q(\eta+1)+2q-2}y]]\\
&=&  [v_{n-2} y\gr{x}{2q-2}(y\gr{x}{2q-1})^{\eta-2}y\gr{x}{2q-2}
 [y\gr{z}{2q(\eta+1)+2q-2}]yx]\\
&\phantom{=}&+  [v_{n-2} y\gr{x}{2q-2}(y\gr{x}{2q-1})^{\eta-2}y\gr{x}{2q-2}x
 [y\gr{z}{2q(\eta+1)+2q-2}]y]\\
&\phantom{=}&+  [v_{n-2} y\gr{x}{2q-2}(y\gr{x}{2q-1})^{\eta-2}y\gr{x}{2q-2}xy
 [y\gr{z}{2q(\eta+1)+2q-2}]]\\
&=& a_1\cdot [v_{n-2} y\gr{x}{2q-2}(y\gr{x}{2q-1})^{\eta-2}y\gr{x}{2q-2} \gr{z}{2q(\eta+1)+2q-1}yx]  \\
&\phantom{=}&+ a _0\cdot [v_{n-2} y\gr{x}{2q-2}(y\gr{x}{2q-1})^{\eta-2}y\gr{x}{2q-2}x \gr{z}{2q(\eta+1)+2q-1}y] \\
&\phantom{=}&+ [v_{n-2} y\gr{x}{2q-2}(y\gr{x}{2q-1})^{\eta-2}y\gr{x}{2q-2}xy \gr{z}{2q(\eta+1)+2q-2}y]\\
&\phantom{=}&+ (b+c_0)\cdot [v_{n-2} y\gr{x}{2q-2}(y\gr{x}{2q-1})^{\eta-2}y\gr{x}{2q-2}xy \gr{z}{2q(\eta+1)+2q-4}yxz]\\
&\phantom{=}&+ (b+c_1)\cdot [v_{n-2} y\gr{x}{2q-2}(y\gr{x}{2q-1})^{\eta-2}y\gr{x}{2q-2}xy \gr{z}{2q(\eta+1)+2q-4}xyz]\\
&=&  [v_{n-2} y\gr{x}{2q-2}(y\gr{x}{2q-1})^{\eta-2}y\gr{x}{2q-2}x
\gr{z}{2q(\eta+1)+2q-1}y] \\
&\phantom{=}&+  [v_{n-2} y\gr{x}{2q-2}(y\gr{x}{2q-1})^{\eta-2}y\gr{x}{2q-2}xy
\gr{z}{2q(\eta+1)+2q-2}y]\\
&\phantom{=}&+  [v_{n-2} y\gr{x}{2q-2}(y\gr{x}{2q-1})^{\eta-2}y\gr{x}{2q-2}xy
\gr{z}{2q(\eta+1)+2q-4}xyz]\\
&=&  [v_{n-2} y\gr{x}{2q-2}(y\gr{x}{2q-1})^{\eta-2}y\gr{x}{2q-2}xy
\gr{z}{2q(\eta+1)+2q-4}xyx]\\
&=&  [\theta_{\frac{n-1}{2}}^\omega x] \ ,
\end{eqnarray*}
where, for $t=0,1$, we have
\begin{eqnarray*}
a_t &=& \binom{2q(\eta+1)+2q-2}{t} + \binom{2q(\eta+1)+2q-2}{2q-1+t} \\
&&\quad + \binom{2q(\eta+1)+2q-2}{2q+2q-2+t}+ \sum_{i=1}^{\eta-1} \binom{2q(\eta+1)+2q-2}{2q(i+1)+2q-2+t}\\
b &=& \binom{2q(\eta+1)+2q-2}{2q-2}  + \binom{2q(\eta+1)+2q-2}{2q+2q-3} + \sum_{i=1}^{\eta-1}  \binom{2q(\eta+1)+2q-2}{2q(i+1)+2q-3} \\
c_t &=& \binom{2q(\eta+1)+2q-2}{2q(\eta+1)+2q-4+t}\ .
\end{eqnarray*}
The identity follows by applying Lucas' Theorem to the above binomial coefficients, since for $0\leq a\leq \eta+1$ and $1\leq b\leq 2q$ we get
\begin{displaymath}
{\binom{2q(\eta+1)+2q-2}{2qa+2q-b}}\equiv {\binom{\eta+1}{a}}{\binom{2q-2}{2q-b}}\equiv {\binom{2^g}{a}}{\binom{2q-2}{2q-b}}\pmod 2\ ,
\end{displaymath}
which does not vanish only when $a=0,\eta+1$ and $b=2,4,2q$.  

\subsection{Expansion of eq. \eqref{theta^ax}}
To prove equation~\eqref{theta^ax}, let $a=h+2-s$, for $1\leq s\leq h$
and expand as follows:
\begin{eqnarray*}
0 &=&  [v_{n-1}y\gr{x}{2q-2}(y\gr{x}{2q-1})^{\eta-2} y\gr{x}{2q-2} [\theta_0^{h+2-s}x]] \\
&=&  [v_{n-1}y\gr{x}{2q-2}(y\gr{x}{2q-1})^{\eta-2} y\gr{x}{2q-2} [y\gr{x}{2q-1}y\gr{x}{2q-2^s-1}yx]] \\
&=&  [v_{n-1}y\gr{x}{2q-2}(y\gr{x}{2q-1})^{\eta-2} y\gr{x}{2q-2}  [y\gr{z}{4q-2^s}x]]\\
&=&  [v_{n-1}y\gr{x}{2q-2}(y\gr{x}{2q-1})^{\eta-2}y\gr{x}{2q-2} [y\gr{z}{4q-2^s}]x]\\
&\phantom{=}&\quad+ [v_{n-1}y\gr{x}{2q-2}(y\gr{x}{2q-1})^{\eta-2}y\gr{x}{2q-2}x[y\gr{z}{4q-2^s}]] \\
&=& 
\left(  \binom{4q-2^s}{1} 
+{\binom{4q-2^s}{2q}}
+{\binom{4q-2^s}{2q+2q-2^s}}
\right)\\
&\phantom{=}&\quad\cdot [v_{n-1}y\gr{x}{2q-2}(y\gr{x}{2q-1})^{\eta-2}y\gr{x}{2q-2} \gr{z}{4q-2^s}yx] \\
&\phantom{=}&+ [v_{n-1}y\gr{x}{2q-2}(y\gr{x}{2q-1})^{\eta-2}y\gr{x}{2q-2}  x  [y\gr{z}{4q-2^s}]] \\
&=& 
\left(  {\binom{4q-2^s}{1}} 
+{\binom{4q-2^s}{2q}}
\right)\cdot [v_{n-1}y\gr{x}{2q-2}(y\gr{x}{2q-1})^{\eta-2}y\gr{x}{2q-2} \gr{z}{4q-2^s}xx] \\
&\phantom{=}&+  
\left(  {\binom{4q-2^s}{0}} 
+{\binom{4q-2^s}{2q-1}}
+{\binom{4q-2^s}{2q+2q-2^s-1}}
\right)\\
&\phantom{=}&\quad\cdot [v_{n-1}y\gr{x}{2q-2}(y\gr{x}{2q-1})^{\eta-1} y\gr{z}{4q-2^s-1}yx]\\
&\phantom{=}&+  
\left(  {\binom{4q-2^s}{0}} 
+{\binom{4q-2^s}{2q-1}}
\right)\cdot
[v_{n-1}y\gr{x}{2q-2}(y\gr{x}{2q-1})^{\eta-1} y\gr{z}{4q-2^s-1}xx]\\
&=&  [v_{n}y\gr{x}{2q-2^s-1}xx] +  [v_{n}y\gr{x}{2q-2^s-1}yx] +  [v_{n}y\gr{x}{2q-2^s-1}xx] \\
&=&  [v_{n}y\gr{x}{2q-2^s-1}yx] \\
&=&  [\theta_1^{h+2-s}x] \ ,
\end{eqnarray*}
since the binomial coefficients which are not immediate can be evaluated by means of Lucas' Theorem as follows:
\begin{eqnarray*}
 {\binom{4q-2^s}{2q}} &\equiv&  {\binom{4q-2^s}{2q+2q-2^s}}\equiv 1 \pmod 2  \ ,\\
 {\binom{4q-2^s}{2q-1}} &\equiv&  {\binom{4q-2^s}{2q+2q-2^s-1}} \equiv 0\pmod 2\ .  
\end{eqnarray*}

\subsection{Expansion of eq. \eqref{Xi}}   
First we expand in the case $n$ even:
\begin{eqnarray*}
0 &=& [ v_s v_s ]\\
&=& [ v_s  [y\gr{z}{2q-1+ds}]]\\
&=& [v_s \gr{z}{2q-1+ds} y] +  
\left( \sum_{l=0}^{s-1} \left( {\binom{2q-1+ds}{dl}} + \sum_{j=0}^{\eta-1} {\binom{2q-1+ds}{2q-1+2qj+dl}} \right) \right.\\ 
&&\quad \left. + {\binom{2q-1+ds}{ds}} \right)
\cdot [v_s \gr{z}{2q+ds}]\\
&=& [v_n y\gr{x}{2q-1}] \ .
\end{eqnarray*}
In fact, the involved binomial coefficients can be rewritten first as
\begin{displaymath}
 \binom{2q-1+ds}{0} + \sum_{l=1}^{s}\left(
 {\binom{2q-1+ds}{dl}}+ \sum_{j=0}^{\eta-1}  {\binom{2q-1+ds}{2q-1+2qj+d(l-1)}}
\right)\ ,
\end{displaymath}
and then, by using binomial properties, as
\begin{displaymath}
1 + \sum_{l=1}^{s}\left( {\binom{2q-1+ds}{dl}} +  {\binom{2q-1+ds}{d(s-l+1)}}
+ \sum_{j=1}^{\eta-1}  {\binom{2q-1+ds}{d(s-l+1)-2qj}}
\right)\ ;
\end{displaymath}
but now,
\begin{displaymath}
\sum_{l=1}^{s}  {\binom{2q-1+ds}{dl}} = \sum_{l=1}^{s} {\binom{2q-1+ds}{d(s-l+1)}}\ ,
\end{displaymath}
and so the remaining coefficient is just
\begin{displaymath}
1 + \sum_{l=1}^{s}\sum_{j=1}^{\eta-1}  {\binom{2q-1+ds}{d(s-l+1)-2qj}}=
1 + \sum_{l=1}^{s}\sum_{j=1}^{\eta-1}  {\binom{2q-1+ds}{dl-2qj}}\ .
\end{displaymath}
Now, by Lucas' Theorem and ($\mathcal{I}$), one gets
\begin{displaymath}
\begin{split}
\sum_{l=1}^{s}\sum_{j=1}^{\eta-1}  {\binom{2q-1+ds}{dl-2qj}} 
&\equiv \sum_{l=1}^{s}\sum_{j=1}^{2^g-2}  {\binom{(2^{g+h+1}-2)s+2^{h+1}-1}{(2^{g+h+1}-2)l- 2^{h+1}j}} \pmod 2\\
&\equiv \sum_{l=1}^{s}\sum_{j=1}^{2^g-2}  {\binom{(2^{g+h+1}-2)s+2^{h+1}-2}{(2^{g+h+1}-2)l- 2^{h+1}j}} \pmod 2\\
&\equiv \sum_{l=1}^{s}\sum_{j=1}^{2^g-2}  {\binom{(2^{g+h}-1)s+2^{h}-1}{(2^{g+h}-1)l- 2^{h}j}} \pmod 2\\
&\equiv \sum_{l=0}^{s-1}\sum_{j=1}^{2^g-2}  {\binom{(2^{g+h}-1)s+2^{h}-1}{(2^{g+h}-1)l+(2^{g+h}- 2^{h}j-1)}} \pmod 2\\
&\equiv 0 \pmod 2\ ,
\end{split}
\end{displaymath}
since $2^{h}-1<2^{g+h}- 2^{h}j-1$ when $j$ ranges between $1$ and $2^g-2$.

Then we expand equation~\eqref{Xi} in the $n$ odd case:
\begin{eqnarray*}
0 &=& [  [v_s y\gr{x}{2q-2}
(y\gr{x}{2q-1})^{\frac{\eta-1}{2}}] 
 [v_s y\gr{x}{2q-2}
(y\gr{x}{2q-1})^{\frac{\eta-1}{2}}] ] \\
&=& [  [v_s y\gr{x}{2q-2} (y\gr{x}{2q-1})^{\frac{\eta-1}{2}}]  
 [y\gr{z}{2q\left( \frac{\eta+1}{2}\right) + 2q-2+ds}]] \\
&=& [v_s y\gr{x}{2q-2} (y\gr{x}{2q-1})^{\frac{\eta-1}{2}} \gr{z}{2q\left( \frac{\eta+1}{2}\right) + 2q-2+ds} y] \\
&\phantom{=}&+  [v_s y\gr{x}{2q-2} (y\gr{x}{2q-1})^{\frac{\eta-1}{2}} \gr{z}{2q\left( \frac{\eta+1}{2}\right) + 2q-1+ds}] \cdot \\
&&\quad \cdot \left( \sum_{j=0}^{\frac{\eta-1}{2}}  \binom{2q\left( \frac{\eta+1}{2}\right) + 2q-2+ds}{2qj} + {\binom{2q\left(\frac{\eta+1}{2}\right) + 2q-2+ds}{2q\left(\frac{\eta-1}{2}\right)+2q-1}} \right. \\
&&\quad\left. + \sum_{l=0}^{s-1} \left( \sum_{j=0}^{\eta-1} {\binom{2q\left(\frac{\eta+1}{2}\right) + 2q-2+ds}{2q\left(\frac{\eta+1}{2}+j\right)+2q-2+dl}} +  {\binom{2q\left(\frac{\eta+1}{2}\right) + 2q-2+ds}{2q\left(\frac{\eta+1}{2}\right)+d(l+1)-1}} \right) \right)\\
&=& [v_s y\gr{x}{2q-2} (y\gr{x}{2q-1})^{\frac{\eta-1}{2}} \gr{z}{2q\left( \frac{\eta+1}{2}\right) + 2q-2+ds} y] \\
&\phantom{=}&+  [v_s y\gr{x}{2q-2} (y\gr{x}{2q-1})^{\frac{\eta-1}{2}} \gr{z}{2q\left( \frac{\eta+1}{2}\right) + 2q-1+ds}] \cdot \\
&&\quad \cdot \left( \sum_{j=0}^{\frac{\eta-1}{2}}  \binom{2q\left( \frac{\eta+1}{2}\right) + 2q-2+ds}{2qj} + \sum_{l=0}^{s-1} \sum_{j=0}^{\eta-1} {\binom{2q\left(\frac{\eta+1}{2}\right) + 2q-2+ds}{2q\left(\frac{\eta+1}{2}+j\right)+2q-2+dl}} \right)\\
&=& [v_s y\gr{x}{2q-2} (y\gr{x}{2q-1})^{\frac{\eta-1}{2}} \gr{z}{2q\left( \frac{\eta+1}{2}\right) + 2q-2+ds} y] \\
&\phantom{=}&+S\cdot  [v_s y\gr{x}{2q-2} (y\gr{x}{2q-1})^{\frac{\eta-1}{2}} \gr{z}{2q\left( \frac{\eta+1}{2}\right) + 2q-1+ds}] \\
&=& [v_n y\gr{x}{2q-1}]\ .
\end{eqnarray*}
The second and the last binomial coefficients vanish since they have an even numerator and an odd denominator.
The sum $S$ of the two remaining terms can be arranged as follows, by using Lucas' Theorem:
\begin{displaymath}
\begin{split}
S &= \sum_{j=0}^{\frac{\eta-1}{2}}  \binom{2q\left( \frac{\eta+1}{2}\right) + 2q-2+ds}{2qj} + 
\sum_{l=0}^{s-1} \sum_{j=0}^{\eta-1} {\binom{2q\left(\frac{\eta+1}{2}\right) + 2q-2+ds}{2q\left(\frac{\eta+1}{2}+j\right)+2q-2+dl}} \\
&= \sum_{j=0}^{2^{g-1}-1} \binom{(2^{g+h+1}-2)s+(2^{g+h}+2^{h+1}-2)}{2^{h+1}j} \\
&\quad+ \sum_{l=0}^{s-1} \sum_{j=0}^{2^g-2} \binom{(2^{g+h+1}-2)s+2^{g+h}+2^{h+1}-2}{(2^{g+h+1}-2)l+2^{g+h}+2^{h+1}(j+1)-2} \\
&= \sum_{j=0}^{2^{g-1}-1} \binom{(2^{g+h}-1)s+2^{g+h-1}+2^{h}-1}{2^{h}j}\\ 
&\quad+\sum_{l=0}^{s-1} \sum_{j=0}^{2^g-2} \binom{(2^{g+h}-1)s+2^{g+h-1}+2^{h}-1}{(2^{g+h}-1)l+2^{g+h-1}+2^{h}(j+1)-1} \\
&= \sum_{j=0}^{2^{g-1}-1} \binom{(2^{g+h}-1)s+2^{g+h-1}+2^{h}-1}{2^{h}j} \\
&\quad+\sum_{l=0}^{s-1} \sum_{j=0}^{2^{g-1}-2} \binom{(2^{g+h}-1)s+2^{g+h-1}+2^{h}-1}{(2^{g+h}-1)l+2^{g+h-1}+2^{h}(j+1)-1}\\ 
&\quad+\sum_{l=0}^{s-1} \sum_{j=2^{g-1}-1}^{2^g-2} \binom{(2^{g+h}-1)s+2^{g+h-1}+2^{h}-1}{(2^{g+h}-1)l+2^{g+h-1}+2^{h}(j+1)-1} \\
&= \sum_{j=0}^{2^{g-1}-1} \binom{(2^{g+h}-1)s+2^{g+h-1}+2^{h}-1}{2^{h}j} \\
&\quad+\sum_{l=0}^{s-1} \sum_{j=0}^{2^{g-1}-2} \binom{(2^{g+h}-1)s+2^{g+h-1}+2^{h}-1}{(2^{g+h}-1)l+2^{g+h-1}+2^{h}(j+1)-1}\\ 
&\quad+ \sum_{l=0}^{s-1} \sum_{j=0}^{2^{g-1}-1} \binom{(2^{g+h}-1)s+2^{g+h-1}+2^{h}-1}{(2^{g+h}-1)(l+1)+2^{h}j} \\
\end{split}
\end{displaymath}
The first term can be grouped into the last one, and then $S$ can be evaluated by using the property ($\mathcal{I}$):
\begin{displaymath}
\begin{split}
S &= \sum_{l=0}^{s-1} \sum_{j=0}^{2^{g-1}-2} \binom{(2^{g+h}-1)s+2^{g+h-1}+2^{h}-1}{(2^{g+h}-1)l+2^{g+h-1}+2^{h}(j+1)-1} \\ 
&\quad+ \sum_{l=0}^{s} \sum_{j=0}^{2^{g-1}-1} \binom{(2^{g+h}-1)s+2^{g+h-1}+2^{h}-1}{(2^{g+h}-1)l+2^{h}j} \\
&= \sum_{l=0}^{s-1} \binom{(2^{g+h}-1)s+2^{g+h-1}+2^{h}-1}{(2^{g+h}-1)l+2^{g+h-1}+2^{h}-1}\\ 
&\quad+ \sum_{l=0}^{s-1} \sum_{j=1}^{2^{g-1}-2} \binom{(2^{g+h}-1)s+2^{g+h-1}+2^{h}-1}{(2^{g+h}-1)l+2^{g+h-1}+2^{h}(j+1)-1}  \\
&\quad+ \sum_{l=0}^{s} \sum_{j=0}^{2^{g-1}-1} \binom{(2^{g+h}-1)s+2^{g+h-1}+2^{h}-1}{(2^{g+h}-1)l+2^{h}j} \\
&= \sum_{l=0}^{s} \binom{(2^{g+h}-1)s+2^{g+h-1}+2^{h}-1}{(2^{g+h}-1)l+2^{g+h-1}+2^{h}-1} - \binom{(2^{g+h}-1)s+2^{g+h-1}+2^{h}-1}{(2^{g+h}-1)s+2^{g+h-1}+2^{h}-1} \\
&\quad +\sum_{l=0}^{s-1} \sum_{j=1}^{2^{g-1}-2} \binom{(2^{g+h}-1)s+2^{g+h-1}+2^{h}-1}{(2^{g+h}-1)l+2^{g+h-1}+2^{h}(j+1)-1} \\ 
&\quad+\sum_{l=0}^{s} \sum_{j=0}^{2^{g-1}-1} \binom{(2^{g+h}-1)s+2^{g+h-1}+2^{h}-1}{(2^{g+h}-1)l+2^{h}j} \\
&\equiv \binom{2^{g+h-1}+2^{h}-1}{2^{g+h-1}+2^{h}-1} - 1 + 0 + \sum_{j=0}^{2^{g-1}-1} \binom{(2^{g+h}-1)s+2^{g+h-1}+2^{h}-1}{(2^{g+h}-1)l+2^{h}j} \pmod 2\\
&\equiv \sum_{j=0}^{2^{g-1}-1} \binom{2^{g+h-1}+2^{h}-1}{2^{h}j} \pmod 2\\
&\equiv 1 + \sum_{j=1}^{2^{g-1}-1} \binom{2^{g+h-1}}{2^{h}j}\binom{2^{h}-1}{0} \pmod 2\\
&\equiv 1\pmod 2\ ,
\end{split}
\end{displaymath}
and thus the thesis.

\subsection{Expansion of eq. \eqref{pag15}}
Then we proceed with~\eqref{pag15}, where $i=\eta-2^{g+h+1-b}$:
\begin{eqnarray*}
0 &=&  [v_{n-1}y\gr{x}{2q-2}(y\gr{x}{2q-1})^{\eta-2}y\gr{x}{2q-2}
 [\theta_0^b x]] \\ 
 &=& [v_{n-1}y\gr{x}{2q-2}(y\gr{x}{2q-1})^{\eta-2}y\gr{x}{2q-2}
 [y\gr{x}{2q-1}y\gr{x}{2q-2}(y\gr{x}{2q-1})^{i}y\gr{x}{2q-2}yx]] \\
&=&  [v_{n-1}y\gr{x}{2q-2}(y\gr{x}{2q-1})^{\eta-2}y\gr{x}{2q-2}
 [y\gr{z}{2q-2+2q(i+2)}x]]\\
&=&  [v_{n-1}y\gr{x}{2q-2}(y\gr{x}{2q-1})^{\eta-2}y\gr{x}{2q-2}
 [y\gr{z}{2q-2+2q(i+2)}]x]\\
&\phantom{=}&+ 
 [v_{n-1}y\gr{x}{2q-2}(y\gr{x}{2q-1})^{\eta-2} y\gr{x}{2q-2}x
 [y\gr{z}{2q-2+2q(i+2)}]]\\
&=& r_1\cdot [v_{n-1}y\gr{x}{2q-2}(y\gr{x}{2q-1})^{\eta-2}y\gr{x}{2q-2} \gr{z}{2q-3+2q(i+2)}yx] \\
&\phantom{=}&+  [v_{n-1}y\gr{x}{2q-2}(y\gr{x}{2q-1})^{\eta-1}\gr{z}{2q-2+2q(i+2)}y]\\
&\phantom{=}&+ r_0\cdot [v_{n-1}y\gr{x}{2q-2}(y\gr{x}{2q-1})^{\eta-1}\gr{z}{2q-1+2q(i+2)}]\\
&=& [v_{n-1}y\gr{x}{2q-2}(y\gr{x}{2q-1})^{\eta-1}\gr{z}{2q-2+2q(i+2)}y]+  [v_{n-1}y\gr{x}{2q-2}(y\gr{x}{2q-1})^{\eta-1}\gr{z}{2q-1+2q(i+2)}]\\
&=& [v_{n-1}y\gr{x}{2q-2}(y\gr{x}{2q-1})^{\eta-1}\gr{z}{2q-3+2q(i+2)}yx]\\
&=& [\theta_n^b x]\ ,
\end{eqnarray*}
where for $t=0,1$
\begin{eqnarray*}
r_t &=& \binom{2q(i+2)+2q-2}{t} + \binom{2q(i+2)+2q-2}{2q-1+t} + \binom{2q(i+2)+2q-2}{2q+2q-2+t} \\
&\phantom{=}&+ \sum_{j=1}^{i} \binom{2q(i+2)+2q-2}{2q(j+1)+2q-2+t} + \binom{2q(i+2)+2q-2}{2q(i+2)+2q-3+t}\ .
\end{eqnarray*}
Since the binomial coefficients whose denominator is odd vanish, 
the non-zero and non-trivial binomial coefficients reduce to
\begin{eqnarray*}
 {\binom{2q(i+2)+2q-2}{2q}} &\equiv&  {\binom{i+2}{1}}=\eta-2^{g+h+1-b}+2
\equiv 1 \pmod 2\ , \\
 {\binom{2q(i+2)+2q-2}{2q(j+1)+2q-2}}  &\equiv&  {\binom{i+2}{j+1}} \pmod 2\ ,
\end{eqnarray*}
and thus the identity follows.

\subsection{Expansion of eq. \eqref{pag16}}
\begin{eqnarray*}
0 &=& [v_{n-1}y\gr{x}{2q-2}(y\gr{x}{2q-1})^{\eta-2}y\gr{x}{2q-1-2^s} [\mu_{0,2} y] ] \\
 &=& [v_{n-1}y\gr{x}{2q-2}(y\gr{x}{2q-1})^{\eta-2}y\gr{x}{2q-1-2^s} [y\gr{x}{2q-1}y\gr{x}{2q-2}y\gr{x}{2q-2}y]] \\
&=& [v_{n-1}y\gr{x}{2q-2}(y\gr{x}{2q-1})^{\eta-2}y\gr{x}{2q-1-2^s} [y\gr{z}{6q-3}y]] \\
&=& [v_{n-1}y\gr{x}{2q-2}(y\gr{x}{2q-1})^{\eta-2}y\gr{x}{2q-1-2^s} [y\gr{z}{6q-3}]y] \\
&=& A\cdot [v_{n-1}y\gr{x}{2q-2}(y\gr{x}{2q-1})^{\eta-2}y\gr{x}{2q-1-2^s}\gr{z}{6q-2}y]\\
&=& [v_{n-1}y\gr{x}{2q-2}(y\gr{x}{2q-1})^{\eta-2}y\gr{x}{2q-1}y\gr{x}{2q-2} y\gr{x}{2q-2}y\gr{x}{2q-2^s-1}y] \\
&=& [v_n y\gr{x}{2q-2}y\gr{x}{k-1}y]\ ,
\end{eqnarray*}
for 
\begin{displaymath}
A = \left( {\binom{6q-3}{2^s}} + {\binom{6q-3}{2q+2^s-1}} + {\binom{6q-3}{4q+2^s-2}} \right)\ .
\end{displaymath}
In fact, the evaluation of the coefficient $A$ can be done as usual by
Lucas' Theorem first
\begin{eqnarray*}
{\binom{6q-3}{2^s}} &\equiv& {\binom{4q}{0}}{\binom{2q-3}{2^s}}\equiv
 {\binom{2q-3}{2^s}} \pmod 2\ ,\\
 {\binom{6q-3}{2q+2^s-1}} &\equiv&
 {\binom{4q}{0}} {\binom{0}{2q}} {\binom{2q-3}{2^s-1}} 
\equiv 1\cdot 0 \cdot  \binom{2q-3}{2^s-1} \equiv 0 \pmod 2\ , \\
 {\binom{6q-3}{4q+2^s-2}} &\equiv&  {\binom{4q}{4q}} {\binom{2q-3}{2^s-2}} \equiv
 {\binom{2q-3}{2^s-2}} \pmod 2\ ,\\
\end{eqnarray*}
and then by using some elementary properties
\begin{eqnarray*}
 {\binom{2q-3}{2^s}} +  {\binom{2q-3}{2^s-2}} &\equiv&  {\binom{2q-3}{2^s}} +
2 {\binom{2q-3}{2^s-1}} +  {\binom{2q-3}{2^s-2}} \pmod 2 \\
&=& \left( {\binom{2q-3}{2^s}}+ {\binom{2q-3}{2^s-1}}+ {\binom{2q-3}{2^s-1}}+ {\binom{2q-3}{2^s-2}}\right)\\
&=& {\binom{2q-2}{2^s}}+ {\binom{2q-2}{2^s-1}}\\
&=&  {\binom{2q-1}{2^s}}\\
&\equiv& 1\pmod 2\ .
\end{eqnarray*}

\subsection{Expansion of eq. \eqref{long2}}
\begin{eqnarray*}
0 &=&  [\mu_{n,i}  [\theta_0^{s+1}x]]] \\
&=& [\mu_{n,i}  [y\gr{x}{2q-1}y\gr{x}{2q-2^{h+2-(s+1)}}x]] \\
&=&  [\mu_{n,i}  [y\gr{z}{4q-2^{h+1-s}}x]]\\
&=&  [\mu_{n,i}  [y\gr{z}{4q-2^{h+1-s}}]x] + [\mu_{n,i} x  [y\gr{z}{4q-2^{h+1-s}}]]\\
&=&  [\mu_{n,i} \gr{z}{4q-2^{h+1-s}+1}x] \cdot\left( {\binom{4q-2^{h+1-s}}{1}} +{\binom{4q-2^{h+1-s}}{2q+1}}\right) \\
&\phantom{=}&+  [\mu_{n,i} x \gr{z}{4q-2^{h+1-s}}y] \\
&\phantom{=}&+  [\mu_{n,i} x \gr{z}{4q-2^{h+1-s}+1}] \cdot\left( {\binom{4q-2^{h+1-s}}{0}} +{\binom{4q-2^{h+1-s}}{2q}}\right) \\
&=& [\mu_{n,i} x \gr{z}{4q-2^{h+1-s}}y] \\
&=& [v_n y\gr{x}{2q-2}(y\gr{x}{2q-1})^i y\gr{x}{k-1}y]\ ,
\end{eqnarray*}
since via Lucas' Theorem the involved binomial coefficients can be
evaluated as follows, where $\alpha, \beta\in\{0,1\}$
\begin{displaymath}
{\binom{4q-2^{h+1-s}}{2q\alpha+\beta}} \equiv
{\binom{1}{\alpha}}{\binom{2q-2^{h+1-s}}{\beta}}\equiv (1-\beta)\pmod 2\ .
\end{displaymath}

\subsection{Expansion of eq. \eqref{lambda}}
Consider the index $0\leq i\leq \eta-3$ and let $\lambda$ be the exponent of the higher
power of two dividing $i+1$, hence $i+1=2^\lambda \cdot r$ where
$\lambda\geq 0$ (with equality when $i$ is even) and $r$ is positive \emph{odd}
integer. Due to the bounds on $i$, we have
\begin{displaymath}
r\leq \left [ \frac{2^g-2}{2^\lambda} \right] \leq  2^{g-\lambda}-1,
\end{displaymath}
where the first becomes an equality if and only if
$i=\eta-2^{g+h+1-b}$ for some $h+2\leq b\leq g+h$, or,
equivalently, $i=2^g-2^\gamma-1$ for some $1\leq \gamma \leq g-1$ so
that $\lambda=\gamma$. 
When we are not in the above case, then
\begin{displaymath}
i+2^\lambda-1
= 2^\lambda\cdot(r+1) -2 
< 2^\lambda \cdot (2^{g-\lambda}-1+1)-2 
\leq 2^g-2^\lambda-2\leq \eta-2\ .
\end{displaymath}
Moreover, if $i$ is even, then $\lambda=0$ and then $i+2^\lambda-1$ is
even, too; otherwise, if $i$ is odd, then $\lambda\geq 1$, and then
$i+2^\lambda-1$ is still even. But $\eta-2^\alpha$ is always even, for
$1\leq \alpha
\leq g-1$, hence $\eta- (i+2^\lambda-1)$ is never $2^\alpha$, for
$1\leq \alpha\leq g-1$. Then $[\mu_{0,i+2^\lambda+1}y]=0$ and we can expand the Jacobi identity
\begin{eqnarray*}
0
&=&[v_{n-1}y\gr{x}{2q-2}(y\gr{x}{2q-1})^{\eta-1-2^\lambda}y\gr{x}{2q-2}
[\mu_{0,i+2^\lambda+1}y]] \\
&=& [v_{n-1}y\gr{x}{2q-2}(y\gr{x}{2q-1})^{\eta-1-2^\lambda}y\gr{x}{2q-2} 
 [y\gr{z}{2q(i+1+2^\lambda)+2q-3}y]] \\
&=&[v_{n-1}y\gr{x}{2q-2}(y\gr{x}{2q-1})^{\eta-1-2^\lambda}y\gr{x}{2q-2}
 [y\gr{z}{2q(i+1+2^\lambda)+2q-3}]y] \\
&=&
 [v_{n-1}y\gr{x}{2q-2}(y\gr{x}{2q-1})^{\eta-1-2^\lambda}y\gr{x}{2q-2}\gr{z}
{2q(i+1+2^\lambda)+2q-3}y] \\
&\qquad&\cdot \left(
  \sum_{j=0}^{2^\lambda-1} {\binom{2q(i+1+2^\lambda)+2q-3}{2qj+1}} + {\binom{2q(i+1+2^\lambda)+2q-3}{2q\cdot 2^\lambda}}\right.\\
&\qquad&\quad \left( +  \sum_{j=0}^{i} {\binom{2q(i+1+2^\lambda)+2q-3}{2q(2^\lambda+j)+2q-1}} \right)\\
&=& [v_n y\gr{x}{2q-2}(y\gr{x}{2q-1})^i y\gr{x}{2q-2}y]\ .
\end{eqnarray*}
As usual, the terms in the last sum vanish since $2q-1>2q-3$, while the remaining
terms can be rewritten together as follows:
\begin{eqnarray*}
\sum_{j=0}^{2^\lambda} {\binom{2q(i+1+2^\lambda)}{2qj}} &\equiv& \sum_{j=0}^{2^\lambda} {\binom{2^\lambda\cdot r+2^\lambda}{j}}\\
&\equiv& 1 + \sum_{j=1}^{2^\lambda-1} {\binom{2^\lambda\cdot (r+1)}{j}} +  {\binom{r+1}{1}}\\ 
&\equiv& 1 + 0 + r +1 \\
&\equiv& 1 \pmod 2\ ,
\end{eqnarray*}
since $r$ is odd. 
\end{document}